\documentclass[11pt]{amsart}
\usepackage[marginratio=1:1]{geometry}
\usepackage{color}
\definecolor{refs}{rgb}{0.7,0,0}
\definecolor{ext}{RGB}{112,112,112}
\definecolor{cite}{RGB}{034,113,179}
\usepackage[colorlinks=true,citecolor=cite,linkcolor=refs,urlcolor=ext,backref=page]{hyperref}

\newtheorem{theorem}{Theorem}[section]
\newtheorem{lemma}[theorem]{Lemma}

\newtheorem{proposition}[theorem]{Proposition}

\theoremstyle{definition}

\newcommand{\mycomment}[1]{}

\title[Conformal geodesics and gauged Killing--Yano forms]{Conformal geodesics and gauged Killing--Yano forms in dimension four}

\author[W.~Kry\'nski]{Wojciech Kry\'nski}
\address{Institute of Mathematics, Polish Academy of Sciences, ul. \'Sniadeckich 8, 00-656 Warszawa, Poland}

\email{krynski@impan.pl}


\newcommand{\C}{\mathbb{C}}

\newcommand{\tr}{\operatorname{tr}}

\newcommand{\spn}{\mathrm{span}}

\newcommand{\cA}{\mathcal{A}}
\newcommand{\cO}{\mathcal{O}}
\newcommand{\CP}{\mathbb{CP}}

\begin{document}
\date{\today}
\maketitle

\begin{abstract}
We study the relation between conformal geodesics and conformal Killing-Yano (CKY) 2-forms on anti-self-dual conformal manifolds of dimension 4. In the flat case, twistor quadrics determine self-dual CKY forms whose zero loci contain conformal geodesics. To extend this construction to curved manifolds, we consider self-dual 2-forms with values in a line bundle and satisfying the gauged CKY equation, coupled to connections with anti-self-dual curvature. In this way we provide a local ambient characterization of non-null conformal geodesics as zero curves of gauged CKY forms, with the curvature vanishing along these curves.
\end{abstract}

\section{Introduction}
A recent paper \cite{Araneda} exploits quadrics in the twistor space of flat four-dimensional space $M$ to construct families of non-flat Einstein metrics. A key ingredient is a correspondence between twistor quadrics and special symplectic forms defined on open subsets of $M$. This correspondence is a flat-space instance of Pontecorvo's twistor characterization of  Kähler metrics with vanishing scalar curvature in anti-self-dual conformal classes \cite{Pontecorvo}. An elementary calculation shows that the symplectic forms are obtained by rescaling self-dual solutions of the conformal Killing--Yano (CKY) equation. This normalization is a manifestation of the construction provided in \cite{DunajskiTod}.

Twistor quadrics also encode conformal geodesics \cite{HollandSparling}, which are special curves defined by a conformally invariant third-order ODE (see, e.g., \cite{BE}). More precisely, each of the two rulings of a non-singular quadric in the twistor space determines a non-null conformal geodesic in $M$. 

The two observations establish a link between self-dual CKY two-forms and conformal geodesics in flat space. Locally, the corresponding conformal geodesics can be recovered as curves that lie in the singular set of the associated normalized symplectic forms. This phenomenon fits into a general scheme where $k$-forms satisfying the CKY equations induce conformal geodesics \cite{DT2,GS,Tod}. In general, however, the CKY equations yield only special conformal geodesics. For instance, related zero-locus theorems for normal CKY solutions, under additional algebraic hypotheses, are established in \cite[Proposition 7.1]{GoverSnellTaghavi} and \cite[Theorem 7.3]{CurryGoverSnell}. The former concerns the simultaneous vanishing of a two-form $\omega$ and its exterior derivative $d\omega$, while the latter treats the zero sets of suitable normal CKY forms in Riemannian signature. Another link between CKY two-forms and conformal geodesics is provided by the first integrals used in \cite[eq.\,(1.3)]{DT2} to study the integrability of the conformal geodesic flow on gravitational instantons.

In this note, we consider general anti-self-dual 4-dimensional manifolds $(M,[g])$. We work locally, using complexified geometry unless a real form is specified. Our objective is to extend the correspondence between conformal geodesics and the self-dual CKY forms. The conformal geodesic equation can be written as (see \cite{BE,DK})
\begin{equation}\label{eq:confgeodesics}
 \nabla_VA-\frac{3g(V,A)}{g(V,V)}A+\frac{3g(A,A)}{2g(V,V)}V-g(V,V)P^\sharp(V)+2P(V,V)V=0,
\end{equation}
where for a given curve $\gamma$ we denote $V=\dot\gamma$ and $A=\nabla_VV$. Moreover, $\nabla$ is the Levi-Civita connection of a (not necessarily flat) metric $g$ and $P=\tfrac{1}{2}Ric(g)-\tfrac{R_g}{12}g$ is the Schouten tensor of $g$. 

As mentioned, in the flat case we get a correspondence between solutions to \eqref{eq:confgeodesics} and self-dual solutions to CKY equation. In the curved case, however, we replace the CKY equation by its generalized gauged version. Let $E\to U$ be a complex line bundle over an open subset $U\subset M$, equipped with a connection $\nabla^E$, and let $\omega\in\Gamma(\Omega^2(U)\otimes E)$ be a 2-form with values in $E$. Denote by $D$ the connection induced by $\nabla^E$ and the Levi-Civita connection $\nabla$ of a representative $g\in [g]$. In a local trivialization of $E$, write $\nabla^E=d+\cA$, where $\cA$ is a 1-form, so that
\[
D_X\omega=\nabla_X\omega+\cA(X)\omega.
\]
We consider the gauged CKY equation
\begin{equation}\label{eq:gCKY}
D_X\omega=\tfrac{1}{3}\iota_X d_{\cA}\omega-\tfrac{1}{3}X^\flat\wedge\delta_{\cA}\omega,  \qquad *\omega=\omega,\\
\end{equation}
for every vector field $X$ on $U$, where in this trivialization
\[
d_{\cA}\omega=d \omega+\cA\wedge \omega,
 \qquad \delta_{\cA}\omega=\delta \omega-\iota_{\cA^\sharp}\omega.
\]
Equation \eqref{eq:gCKY} was studied in \cite[eqs. $(16')$ and $(17)$]{BCK}. 
Note that the earlier non-null self-dual case goes back to \cite{DR}, as explained in \cite[Section V]{BCK}. For our construction, we will additionally require the curvature of $\nabla^E$ to be anti-self-dual (i.e. satisfying the instanton equation \cite{Donaldson})
\begin{equation}\label{eq:ASD}
f=d\cA,\qquad *f=-f.
\end{equation}
Our main result is as follows
\begin{theorem}\label{thm:main}
Let $(E,\nabla^E,\omega)$ satisfies \eqref{eq:gCKY} and \eqref{eq:ASD} on an anti-self-dual manifold $(M,[g])$. Fix a local trivialization of $E$ with the corresponding $\cA$ defining $\nabla^E$ and assume that $\xi=\delta_\cA\omega$ and $K=\xi^\sharp$ is a non-null vector field with respect to $g$. Then every regular curve component $\gamma$ of $\{\omega=0\}$ is a conformal geodesic satisfying $\dot\gamma\in\spn\{K\}$ provided that $f(x)=0$ for all $x\in\gamma$.
Conversely, any sufficently small segment of a not-null conformal geodesic can be obtained in this way.
\end{theorem}

Equation \eqref{eq:gCKY} is conformally invariant. More precisely, if $\omega$ satisfies \eqref{eq:gCKY} with respect to $g$ and a fixed connection $\nabla^E$, then $\tilde\omega=e^{3\varphi}\omega$ satisfies the same equation with respect to $\tilde g=e^{2\varphi}g$ and the same connection $\nabla^E$ (see \cite{Semmelmann}).

In the next section, Section \ref{sec:confgeodesics}, we provide basic facts concerning conformal geodesics. Then, in Section \ref{sec:flat} we briefly recall the flat case, where the standard CKY equation is sufficient for the determination of the conformal geodesics. In Section \ref{sec:CKY}, we prove Theorem \ref{thm:CKY} which describe properties of the CKY equation in the curved case. Note that even in this case our assumptions differ from the assumptions in \cite{GoverSnellTaghavi} and \cite{CurryGoverSnell}. 

Finally, in Section \ref{sec:gCKY}, we provide the proof of Theorem \ref{thm:main}. On the one hand, we provide explicit calculations showing that the solutions to \eqref{eq:gCKY} imply solutions to the conformal geodesic equations. This part is the content of Theorem \ref{thm:gCKY}. On the other hand, we use twistor methods (in particular the Ward correspondence) to construct solutions to \eqref{eq:gCKY} from a given conformal geodesic. We do not claim uniqueness: there may be more than one $\omega$ for a given conformal geodesic. This is not the case for the standard CKY in the flat case (up to the conformal rescaling). Indeed, the self-dual solutions to CKY form a 10-parameter family of solutions, and up to a conformal constant, they are in a one-to-one correspondence with the 9-parameter family of the conformal geodesics.


\subsection*{Acknowledgments.}
This work was partially supported by the Simons Foundation grant (award no.\,SFI-MPS-T-Institutes-00010825) and from State Treasury funds as part of a task commissioned by the Minister of Science and Higher Education under the project “Organization of the Simons Semesters at the Banach Center - New Energies in 2026-2028” (agreement no. MNiSW/2025/DAP/491).

AI disclosure: ChatGPT was used to assist with literature searches concerning twistor methods, tensor calculations in the proof of Theorem~\ref{thm:CKY}, and language verification of the paper.

\section{Preliminaries on conformal geodesics}\label{sec:confgeodesics}
As mentioned in the Introduction, conformal geodesics on a conformal manifold $(M,[g])$ are defined by the conformally
invariant third-order ODE \eqref{eq:confgeodesics}. They have a natural projective class of parameters: any two distinguished parameters are related by a M{\"o}bius transformation, corresponding to the group action of $\mathrm{GL}(2)$. This follows from the fact that \eqref{eq:confgeodesics} restricted to a curve is the Schwarzian derivative (see \cite{BE,K}; for more information on conformal geodesics see also \cite{DunajskiTod, GoverSnellTaghavi, Lenka, Kruglikov}).

In what follows, however, it will be more convenient to regard conformal geodesics as unparametrized curves. Let $\gamma$ be a non-null curve in $M$ with an arbitrary parameter, and write $V=\dot\gamma$ and $A=\nabla_VV$. We consider the equation
\begin{equation}\label{eq:unpar-confgeodesic}
 V\wedge\left(\nabla_VA-\frac{3g(V,A)}{g(V,V)}A -g(V,V)P^\sharp(V)\right)=0.
\end{equation}

Equation \eqref{eq:unpar-confgeodesic} can be easily rewritten for curves with unit tangent vector $v$ (i.e. $g(v,v)=1$).  Denoting $a=\nabla_v v$, we have $g(a,v)=0$ and equation \eqref{eq:unpar-confgeodesic} becomes
\begin{equation}\label{eq:circle}
\nabla_v a= P^\sharp(v)-\bigl(P(v,v)+g(a,a)\bigr)v.
\end{equation}
In this viewpoint, the parametrization of a curve is dependent on the choice of $g$ in the conformal class. It is only used for calculations and does not have any invariant meaning. However, \eqref{eq:circle} has a very important consequence. Namely, it follows that a metric geodesic for fixed $g$ (i.e. when $a=0$ is taken in \eqref{eq:circle}) is a conformal geodesic exactly when $P^\sharp(v)$ is parallel to $v$. Moreover, one can invert the reasoning to have the following statement.
\begin{proposition}[\cite{BE}, Propositoin 3.3]
A curve is a conformal geodesic if and only if there exists a metric in the conformal class such that the curve is a metric geodesic and its velocity $v$ is parallel to $P^\sharp(v)$.
\end{proposition}
In fact, \cite[Proposition~3.3]{BE} is a stronger result: one can choose a representative metric $g$ for which the curve is a geodesic and $P^\sharp(v)=0$. We will only need the weaker statement formulated above.

\section{Flat case}\label{sec:flat}
In this section we shall briefly recall standard constructions in the flat case which serve as a motivation for further parts of the paper. Let $M\simeq\C^4$ be equipped with a flat metric. Its twistor space $T$ is $\CP^3\setminus\CP^1$, whose points parametrize one of the two families of totally null affine two-planes in $M$. Each point $x\in M$ corresponds to a projective line $L_x$ in twistor space. We choose homogeneous coordinates so that
\[
L_x=\{[as+bt:cs+dt:s:t] \ |\ [s:t]\in\CP^1\}.
\]
where the coefficients $(a,b,c,d)$ can be taken as coordinates on $M$, with $x=(a,b,c,d)$. In these coordinates the metric is given by
\[
g=2(d a\odot d d-d b\odot d c),
\]
where $\odot$ is the symmetric tensor product: $\alpha\odot\beta=\tfrac{1}{2}(\alpha\otimes\beta+\beta\otimes\alpha)$. Then, the standard CKY equation takes the form
\begin{equation}\label{eq:CKY}
\nabla_a\omega_{bc}=\tfrac13(d\omega)_{abc}+g_{ab}\xi_c-g_{ac}\xi_b,
\end{equation}
where
\begin{equation}
\xi_c=\tfrac{1}{3}\sum_d\nabla^d\omega_{dc},
\end{equation}
and $\nabla$ is the Levi-Civita connection. Notice also that $(\delta \omega)_b=-\sum_a\nabla^a\omega_{ab}$.

Now, let $Q$ be a symmetric matrix, $Q=(q_{ij})_{i,j=0}^3$. In the fixed affine chart on the twistor space it defines a quadratic form $\chi_Q$ defined by formula 
\[
 \chi_Q(Z)=Z^TQZ,
\]
where $Z=(Z^0,\ldots,Z^3)\in\C^4$. Its restriction to the line $L_x$ takes the form $\chi_Q|_{L_x}=Cs^2+2Bst+At^2$, where explicitly
\begin{equation}\label{eq:ABC}
\begin{aligned}
 &C=q_{00}a^2+2q_{01}ac+q_{11}c^2+2q_{02}a+2q_{12}c+q_{22},\\
 &B=q_{00}ab+q_{01}(ad+bc)+q_{11}cd+q_{02}b+q_{03}a+q_{12}d+q_{13}c+q_{23},\\
 &A=q_{00}b^2+2q_{01}bd+q_{11}d^2+2q_{03}b+2q_{13}d+q_{33}.
\end{aligned}
\end{equation}
We shall show that 
\begin{equation}\label{eq:flatform}
\omega=A da\wedge dc-B(d a\wedge d d+d b\wedge d c)+Cd b\wedge d d
\end{equation}
solves \eqref{eq:CKY}. 

Indeed, equation \eqref{eq:CKY} for self dual 2-forms becomes
\begin{equation}\label{eq:flatpde}
\begin{gathered}
A_a=A_c=0,\quad C_b=C_d=0,\\
2B_a=A_b,\quad 2B_c=A_d,\quad
2B_b=C_a,\quad 2B_d=C_c,
\end{gathered}
\end{equation}
and one checks that it is satisfied by \eqref{eq:ABC}. Conversely, equations \eqref{eq:flatpde} imply that
$A_{bb}=C_{aa}$, $A_{bd}=C_{ac}$ and $A_{dd}=C_{cc}$ and they are all constant. Integration gives exactly \eqref{eq:ABC}. Notice that \eqref{eq:flatform} with arbitrary $A,B,C$ is a formula for a general self-dual 2-form with respect to metric $g$ (because $\Omega_1=da\wedge dc$, $\Omega_2=(d a\wedge d d+d b\wedge d c)$ and $\Omega_3=d b\wedge d d$ form a basis of the self-dual 2-forms for $g$). Thus every local self-dual CKY 2-form comes from a unique symmetric $Q$. As remarked in the Introduction, there are 10 linear parameters, or 9 after projectivization.

\subsection*{Symplectic normalization and conformal geodesics}
Let
\[
\Delta=B^2-AC=-\tfrac{1}{4}\sum_{a,b}\omega_{ab}\omega^{ab}.
\]
Then
\[
\sigma=\Delta^{-1}\omega
\]
is the symplectic form from \cite{Araneda} (see also \cite{AASS}). It is well defined on the subset $\Delta\neq 0$. Notice that we can recover $\omega$ from $\sigma$ by 
\[
 \omega=\kappa^{-3/4}\omega,
\]
where $\kappa=-\tfrac{1}{4}\sum_{a,b}\sigma_{ab}\sigma^{ab}$. Hence, it follows that the symplectic forms are characterized intrinsically by requiring that the associated $\omega$ satisfies the CKY equation. Moreover, under $Q\mapsto\lambda Q$, one has $\sigma\mapsto\lambda^{-2}\sigma$.

Let $S_Q=\{\chi_Q=0\}$. Then we have that $L_x\subset S_Q$ if and only if $A=B=C=0$, or equivalently $\omega(x)=0$. On the other hand $\{\Delta=0\}$ consists of points $x$ such that $L_x$ is tangent to $S_Q$. Note that the 2-form $\omega$ is smooth and only its symplectic normalization $\sigma$ has singularities. The following result follows from \cite[Section 1.6]{HollandSparling}.

\begin{proposition}\label{prop:flatcircles}
For a complex quadric $S_Q\subset\CP^3$, the zero locus of the corresponding 2-form $\omega$ consists of two non-null conformal geodesics, one for each ruling.
\end{proposition}
Let us point out that on a compatible real slice, nonempty real loci give real conformal geodesics. In particular, with respect to the standard Euclidean flat coordinates these are ordinary circles (or straight lines).

\section{CKY forms on arbitrary anti-self-dual manifolds}\label{sec:CKY}
Now we shall show that Proposition \ref{prop:flatcircles} partially extends to the curved case. However, one only has the result in one direction.
We will need the following lemma.
\begin{lemma}\label{lemma:ASDforms}
Let $\alpha$ be an anti-self-dual 2-form on $M$ and let $X$ be a non-null vector field. Then
\[
\alpha=\frac{1}{g(X,X)}\left(X^\flat\wedge\iota_X\alpha-*\bigl(X^\flat\wedge\iota_X\alpha\bigr)\right).
\]
In particular, $\iota_X\alpha=0$ implies $\alpha=0$.
\end{lemma}
\begin{proof}
Let us complement $X^\flat$ with $\beta_i$, $i=1,2,3$, to an orthonormal co-frame on $M$ and denote $\mathcal{B}=\spn\{\beta_1,\beta_2,\beta_3\}$. We can decompose $\alpha$ in the following way $\alpha=X^\flat\wedge\beta+\eta$, where $\beta\in\mathcal{B}$ and $\eta\in\mathcal{B}^{\wedge 2}$. Note that necessarily $\beta=g(X,X)^{-1}\iota_X\alpha$. Additionally $*$ interchanges $X^\flat\wedge \mathcal{B}$ and $\mathcal{B}\wedge\mathcal{B}$. Thus, the assumption $*\alpha=-\alpha$ implies
\[
\eta=-*(X^\flat\wedge\beta)
\]
and this gives the required formula.
\end{proof}

\begin{theorem}\label{thm:CKY}
Let $\omega$ be a self-dual CKY 2-form on an anti-self-dual manifold $(M,[g])$. Let $\xi=-\tfrac{1}{3}\delta\omega$ and $K=\xi^\sharp$. Then, every regular curve component $\gamma$ of $\{\omega=0\}$, such that $K$ is non-null along $\gamma$, is a conformal geodesic, satisfying $\dot\gamma\in\spn\{K\}$.
\end{theorem}

\begin{proof}
The CKY equation reads
\[
\nabla_X\omega=\frac{1}{3}\iota_X d\omega-\frac{1}{3}X^\flat\wedge\delta\omega
\]
for any vector field $X$. In case of self-dual $\omega$ (i.e. $\omega=*\omega$) it reduces to
\begin{equation}\label{eq:CKYreduced}
\nabla_X\omega=X^\flat\wedge\xi+*(X^\flat\wedge\xi).
\end{equation}
Denote the right hand side of equation \eqref{eq:CKYreduced} by $E_X(\xi)$, i.e. $E(\xi)$ can be considered as a map that sends $X$ to the 2-form $E_X(\xi)$. By assumption $g(K,K)\neq 0$, and it follows that $X$ is in the kernel of the map $X\mapsto E_X(\xi)$ if and only if $X$ is proportional to $K$. Then, it follows from the implicit function theorem that the zero set of $\omega$ is locally a curve $\gamma$ such that $\dot\gamma\in\spn\{K\}$. Indeed, as a self-dual form $\omega$ has 3-components in a basis of self-dual forms: $\omega=A\Omega_1+B\Omega_2+C\Omega_3$ for some locally non-vanishing self-dual $\Omega_i$, $i=1,2,3$. The zero set $\{\omega=0\}$ is locally described as $A=B=C=0$. Note that the differential of the functions $A,B,C$ on this set can be computed as $\nabla\omega$. By earlier considerations this differential has rank 3, since it equals to the mapping $E(\xi)$.

Now, we shall prove that metric $g$ can be rescaled such that for the new metric $\nabla_KK=0$. Equivalently we can consider $\alpha(X,Y)=(\nabla_X\xi)(Y)$ and show that $\alpha(K,.)=0$ along $\gamma$.

Firstly, we will provide some algebraic properties of $\alpha$ along $\gamma$. Denote $\alpha_X=\alpha(X,\cdot)$. Then, commuting covariant derivatives of equation \eqref{eq:CKYreduced} and using $\omega=0$ gives
\begin{equation}\label{eq:eq1}
E_Y(\alpha_X)-E_X(\alpha_Y)=0
\end{equation}
along $\gamma$. Indeed, from the left hand side of \eqref{eq:CKYreduced} we get $\nabla_Y\nabla_X\omega-\nabla_X\nabla_Y\omega=\nabla_{[Y,X]}\omega+R(X,Y)(\omega)$ and the component with the curvature $R(X,Y)$ vanishes due to $\omega=0$. On the other hand, from the right hand side, we get $E_Y(\alpha_X)-E_X(\alpha_Y)+E_{[Y,X]}(\xi)$. The components involving $[Y,X]$ cancel out by the definition of $E$ and we arrive to \eqref{eq:eq1}. 

Now, chose a local orthonormal frame $(e^0,\ldots,e^3)$ on $M$. Taking $Y=e^i$ in \eqref{eq:eq1}, evaluating on $e^i$ and summing over $i$ we get
\begin{equation}\label{eq:eq2}
0=\sum_i \bigl(E_{e_i}(\alpha_X)-E_X(\alpha_{e_i})(e_i)\bigr)=2\alpha_X+(\tr_g\alpha)X^\flat
   +2\iota_X(*\alpha^{skew}),
\end{equation}
where $\alpha^{skew}$ is the skew-symmetric part of $\alpha$. Similarly, taking also $X=e_i$ in \eqref{eq:eq2}, applying the resulting 1-forms to $e_i$ and summing over $i$ we get $6\tr_g\alpha=0$. Then \eqref{eq:eq2} gives that the symmetric part of $\alpha$ vanishes, and the remaining skew-symmetric part is anti-self dual, i.e.
\[
\alpha(X,Y)=-\alpha(Y,X),\qquad *\alpha=-\alpha
\]
which hold along $\gamma$.

Let $a=\nabla_KK$, then since $\xi=K^\flat$, we have $a^\flat=\alpha_K$ and $g(a,K)=\alpha(K,K)=0$, becasue $\alpha$ is skew-symmetric along $\gamma$. In particular, $K(g(K,K))=2g(a,K)=0$, so $K$ has constant length along $\gamma$.

We define $\hat g=e^{2\varphi}g$, where $\varphi$ is a function satisfying
\begin{equation}\label{eq:phi}
\varphi|_\gamma=0,\quad d\varphi|_\gamma=g(K,K)^{-1}a^\flat.
\end{equation}
Such a function exists because $a^\flat(K)=g(a,K)=0$, so the condition $K(\varphi)=0$ is consistent with $d\varphi(K)=0$. The transformation rule gives
\begin{equation}\label{eq:trans}
\hat\nabla_KK =\nabla_KK+2d\varphi(K)K-g(K,K)(d\varphi)^\sharp=0
\end{equation}
which holds on $\gamma$, because of \eqref{eq:phi}. Then, replacing $\omega$ by $\hat\omega=e^{3\varphi}\omega$, the associated one-form $\hat\xi$ takes the form
\[
\hat\xi=e^\varphi\bigl(\xi+\iota_{(d\varphi)^\sharp}\omega\bigr).
\]
We get that $\hat\xi=\xi$ and $\hat K=K$ on $\gamma$. Consequently we get from \eqref{eq:trans} that  $\hat\alpha(\hat K,\cdot)=(\hat\nabla_{\hat K}\hat K)^\flat=0$. This completes the argument that $\gamma$ is a geodesic for $\hat g$. Notice also that the preceding argument shows that $\hat\alpha$ is anti-self-dual on $\gamma$.

We will now omit hats and prove that $K$ is parallel to $P^\sharp(K)$. Notice that since $\alpha$ is anti-self-dual along $\gamma$, the general formula from Lemma \ref{lemma:ASDforms} implies that $\alpha=0$, because $K$ annihilates $\alpha$ on the right hand side.

As before, we shall denote by $R(X,Y)$ the curvature tensor of $g$. Depending on the context it can act on any tensor on $M$. Let $H_X=\nabla_X\alpha$, viewed as a bilinear form. Then, since $\alpha=0$ along $\gamma$ and $K$ is tangent to $\gamma$, we get $H_K=0$. Commuting derivatives of $\xi$ gives
\[
H_X(Y,Z)-H_Y(X,Z)=(R(X,Y)\xi)(Z).
\]
From this, taking $Y=Z=K$, we obtain $H_X(K,K)=0$. Moreover (as already used before) the following holds
\[
R(Y,Z)\omega=E_Z(\alpha_Y)-E_Y(\alpha_Z).
\]
Differentiating this with $X$ and then restricting to $\gamma$ (i.e. to points where $\omega=0$) gives
\begin{equation}\label{eq:eq3}
E_Z\bigl(H_X(Y,\cdot)\bigr)-E_Y\bigl(H_X(Z,\cdot)\bigr)=R(Y,Z)E_X(\xi).
\end{equation}
After a constant rescaling of $\omega$, we can assume $g(K,K)=1$ and that $(K,X,Y,Z)$ is an orthonormal frame.
Then
\[
E_X(\xi)=-(K^\flat\wedge X^\flat+Y^\flat\wedge Z^\flat).
\]
Now, we apply the preceding identity applied to the triples $(X,K,Y)$, $(X,K,Z)$ and $(X,Y,Z)$, and then evaluate the resulting two-forms on $(K,Y)$, $(K,Z)$ and $(Y,Z)$, respectively. Using $H_X(K,K)=0$, we obtain
\[
\begin{aligned}
&-H_X(Y,Y)=R(K,Y,K,Z)- R(K,Y,Y,X),\\
&-H_X(Z,Z)=- R(K,Z,K,Y)-R(K,Z,Z,X),\\
&H_X(Y,Y)+H_X(Z,Z)=0.
\end{aligned}
\]
Finally, from the symmetry $R(K,Y,K,Z)=R(K,Z,K,Y)$ it follows, by adding the above expressions, that
\begin{equation}\label{eq:ricci}
Ric(K,X)=-R(K,Y,Y,X)-R(K,Z,Z,X)=0.
\end{equation}
Note that $X$ was an arbitrary unit vector perpendicular to $K$. Hence, we deduce that $P^\sharp(K)$ is proportional to $K$. Indeed, by \eqref{eq:ricci}, we have $P(K)=\tfrac{1}{2}Ric(K,\cdot) -\tfrac{R_q}{12}K^\flat\simeq K^\flat$. Together with $\nabla_KK=0$, this proves that $K$ satisfies the unparametrized conformal-geodesic equation, which completes the proof.
\end{proof}

\section{Gauged CKY forms}\label{sec:gCKY}
We shall prove here our main Theorem \ref{thm:main}. It splits into two parts. Firstly, we will prove a generalization of Theorem \ref{thm:CKY} to the gauged version of CKY equation. Then we shall apply the twistor methods, in particular the Ward correspondence, to construct a the CKY data from a conformal geodesic.

\begin{theorem}\label{thm:gCKY}
Let $(E,\nabla^E,\omega)$ satisfies \eqref{eq:gCKY} and \eqref{eq:ASD} on an anti-self-dual manifold $(M,[g])$.  Fix a local trivialization of $E$ with the corresponding $\cA$ defining $\nabla^E$ and assume that $\xi=-\tfrac{1}{3}\delta_\cA\omega$ and $K=\xi^\sharp$ is a non-null vector field with respect to $g$. Moreover, let $p\in M$ be such that $\omega_p= 0$. Then in a neighborhood of $p$, the zero set of $\omega$ is a non-null curve $\gamma$, with $\dot\gamma\in\spn\{K\}$. Moreover, $\gamma$ is a conformal geodesic if and only if $f(x)=0$ for all $x\in\gamma$ where $f$ is defined by \eqref{eq:ASD}.
\end{theorem}
\begin{proof}
Let us write \eqref{eq:gCKY} as
\[
D_X\omega=E_X(\xi)
\]
where $E_X(\xi)=X^\flat\wedge\xi+*(X^\flat\wedge\xi)$. We shall also denote $\alpha_X=D_X\xi$. With this notation the proof follows the line of the proof of Theorem \ref{thm:CKY}.

Firstly, at zero set of $\omega$, the same rank argument as in Theorem \ref{thm:CKY} gives a smooth curve $\gamma$ with $\dot\gamma=\spn\{K\}$. Then, the same contraction of $E_X(\alpha)$, computed in a local orthonormal frame as in \eqref{eq:eq2}, proves that $\alpha$ is skew-symmetric and anti-self-dual along $\gamma$. The difference from the ordinary CKY case is that the curvature is of the form $R+f$. Nevertheless, for $x\in\{\omega=0\}$ the additional term $f(X,Y)\omega$ also vanishes. 

Next, we notice that without loss of generality we can assume that $\cA(K)=0$ along $\gamma$, by choosing a local section $e$ of $E$ such that $\nabla^E_K e=0$. Then $\alpha_K=(\nabla_KK)^\flat$, and the conformal change used in the proof of Theorem \ref{thm:CKY} applies without modification. Omitting hats, we get along $\gamma$ the following formulae
\[
\nabla_KK=0,\qquad \alpha=0,\qquad \cA(K)=0.
\]

Consequently, it suffices to prove that $P^\sharp(K)$ is parallel to $K$. Similarly, to the previous case we set $H_X=D_X\alpha$. Thus, in particular, $H_K=0$. However, unlike in Theorem \ref{thm:CKY}, commuting derivatives of $\xi$ gives
\[
H_X(Y,Z)-H_Y(X,Z)=(R(X,Y)\xi)(Z)+f(X,Y)\xi(Z).
\]
and from this we obtain
\[
H_X(K,K)=f(X,K)=-f(K,X).
\]
Moreover, formula \eqref{eq:eq3} is replaced by
\[
E_Z\bigl(H_X(Y,\cdot)\bigr) -E_Y\bigl(H_X(Z,\cdot)\bigr)=\bigl(R(Y,Z)+f(Y,Z) \mathrm{id}\bigr)E_X(\xi).
\]
However, the same triples of vectors and the same evaluations as in the proof of Theorem~\ref{thm:CKY} give
\[
\begin{aligned}
&f(K,X)-H_X(Y,Y)=R(K,Y,K,Z)-R(K,Y,Y,X),\\
&f(K,X)-H_X(Z,Z)=-R(K,Z,K,Y)-R(K,Z,Z,X),\\
&H_X(Y,Y)+H_X(Z,Z)=f(Y,Z).
\end{aligned}
\]
By our assumption \eqref{eq:ASD}: $*f=-f$. It implies, by expressing it in terms of the co-frame dual to $(K,X,Y,Z)$, that $f(Y,Z)=-f(K,X)$. Then, the above equations imply
\[
Ric(K,X)=2f(K,X)-f(Y,Z)=3f(K,X).
\]
Consequently, we get $P(K,X)=\frac{3}{2}f(K,X)$ which hold for all $X$ perpendicular to $K$. It follows that $K$ is proportional to $P^\sharp(K)$ if and only if $f(K,\cdot)=0$. But, since $f$ is anti-self-dual and $K$ is non-null, we can again exploit Lemma \ref{lemma:ASDforms} and get that $\iota_Kf=0$ is equivalent to $f=0$ along $\gamma$.
\end{proof}

Now we proceed to the construction of a solution $(E,D,\omega)$ to \eqref{eq:gCKY} and \eqref{eq:ASD} from a given conformal geodesic. We will start with a lemma. Recall that $T$ denotes the twistor space of $M$, which is a 3-dimensional complex manifold. Its canonical bundle $K_T$ is the bundle of holomorphic 3-forms. Moreover, for a given smooth holomorphic surface $S\subset T$, $\cO_T(S)$ denotes the associated holomorphic line bundle, equipped with a canonical holomorphic section vanishing precisely on $S$, with multiplicity one.  We also choose a holomorphic square root $K_T^{1/2}$ of $K_T$, so that $K_T^{1/2}\otimes K_T^{1/2}\simeq K_T$. Such a square root exists on a sufficiently small neighbourhood of any twistor line $L_x\simeq\CP^1$. Finally, for a holomorphic submanifold $Y\subset Z$, we write $N_{Y/Z}$ for the normal bundle of $Y$ in $Z$.

\begin{lemma}\label{lemma:l1}
Let $(M,[g])$ be a complex analytic ASD conformal 4-dimensional manifold, and let $\gamma$ be a  holomorphic curve in $M$ with non-null tangents. Locally near any point of $\gamma$, there exist a holomorphic line bundle $E$, a holomorphic connection $\nabla^E$ with anti-self-dual curvature, and an $E$-valued self-dual two-form $\omega$ satisfying the gauged CKY equation \eqref{eq:gCKY} and \eqref{eq:ASD}, such that $\omega$ vanishes along $\gamma$ and $\delta_{\cA}\omega$ is non-null.
\end{lemma}\begin{proof}
We shall use the rank-one case of the Ward correspondence (see e.g. \cite{Donaldson} or \cite{CP, Maciek-book}). To apply it, we will first construct a smooth holomorphic surface in $T$ and then use it to define a holomorphic line bundle on a suitable open subset of $T$, whose restriction to each twistor line corresponding to a point in a neighbourhood of $\gamma$ is holomorphically trivial. Then we will use \cite{DunajskiTod} to prove that the CKY equation is satisfied. 

Let $x$ be a point on $\gamma$. Since we have the identification $T_xM\simeq H^0(L_x,N_{L_x/T})$ a tangent vector $\dot\gamma$ can be identified with a section of the normal bundle of $L_x$. By the description of the null cone in \cite{Donaldson}, non-null vectors correspond to the nowhere-vanishing sections of $N_{L_x/T}$. Thus, after shrinking $\gamma$ around one of its points, $S_\gamma=\bigcup_{x\in\gamma}L_x$ is a smooth embedded holomorphic surface. Indeed, the differential of the map $\gamma\times\CP^1\to T$ sending a point $(x,z)$ to a point on the twistor line $L_x$ has rank two. Indeed, one direction is tangent to the twistor line, while the direction along $\gamma$ has a nonzero normal component given by the section induced by $\dot\gamma$. Hence the map is an immersion.

Now, consider $L_x\subset S_\gamma\subset T$ and the associated normal bundle exact sequence: $0\to N_{L_x/S_\gamma}\to N_{L_x/T}\to N_{S_\gamma/T}|_{L_x}\to 0$. Since $L_x$ has trivial normal bundle in $S_\gamma$, and $N_{L_x/T}\simeq \cO(1)^{\oplus2}$, we have the following normal bundle exact sequence
\[
0\to \cO\to \cO(1)\oplus \cO(1)\to \cO_T(S_\gamma)|_{L_x}\to 0.
\]
Consequently $\cO\otimes  \cO_T(S_\gamma)|_{L_x}\simeq\det\bigl(\cO(1)\oplus\cO(1)\bigr)\simeq \cO(2)$.
Hence $O_T(S_\gamma)|_{L_x}\simeq\cO(2)$. Moreover, since $K_{L_x}\simeq \cO(-2)$ and $N_{L_x/T}\simeq \cO(1)^{\oplus2}$, we have
\[
K_T|_{L_x} \simeq K_{L_x}\otimes(\det N_{L_x/T})^{-1} \simeq \cO(-4).
\]
and then
\[
\mathcal{L}_\gamma=\cO_T(S_\gamma)\otimes K_T^{1/2}
\]
has degree zero on lines $L_x$ (and on all twistor lines corresponding to points sufficiently close to $x$). Its restrictions are therefore holomorphically trivial. The Ward correspondence produces a line bundle $E$ over $M$ (on a neighborhood of $\gamma$) with a connection $\nabla^E$ whose curvature is anti-self-dual, i.e.\,satisfies \eqref{eq:ASD}.

Now, since we have $O_T(S_\gamma)|_{L_x}\simeq\cO(2)$, it follows that the canonical section of $\cO_T(S_\gamma)$ restricts to a quadratic section on each twistor line (exactly as in the flat case, when the canonical section is given by the quadric $Q$). Moreover, $\cO_T(S_\gamma)=K_T^{-1/2}\otimes\mathcal{L}_\gamma$ and, after choosing $g\in[g]$, these restrictions define an $E$-valued self-dual 2-form $\omega$. Then, one can adapt the reasoning of \cite[Section~5]{DunajskiTod} for the standard CKY equation to the gauged case and conclude that $(D_X\omega)(Y,Z)=0$ whenever $X,Y,Z$ lie in one $\alpha$-plane. For self-dual $\omega$, this condition is equivalent to the (gauged) CKY equation.

Finally, $\omega_x=0$ for every $x\in\gamma$, since $L_x\subset S_\gamma$ and the canonical section defining $S_\gamma$ vanishes identically on $L_x$. Hence, it remains to show that $\xi=-\tfrac{1}{3}\delta_{\mathcal A}\omega$ is non-null. Let $x\in\gamma$ and $u\in T_xM$. Then by the gauged CKY equation we have
\begin{equation}\label{eq:linear}
(D_u\omega)_x=u^\flat\wedge\xi_x+*(u^\flat\wedge\xi_x).
\end{equation}
In particular, \eqref{eq:linear} implies that $\xi$ completely determines $D_u\omega$ along $\gamma$. Take $v\in T_x\gamma$. Then $D_v\omega=0$, so by contraction of \eqref{eq:linear} with $v$, we obtain
\[
g(v,v)\xi_x-\xi_x(v)v^\flat=0.
\]
It follows that $\xi_x=c v^\flat$ for some $c$. In order to show that $c\neq0$, we consider the
linearization of the section $\omega$, namely the map $u\mapsto(D_u\omega)_x$, at points $x$ where $\omega_x=0$. Since $\omega_x$ is represented by the restriction of the defining section of $S_\gamma$ to $L_x$, its linearization at a zero is induced by the normal bundle map
\[
N_{L_x/T}\to N_{S_\gamma/T}|_{L_x}.
\]
Indeed, differentiating the defining section along a deformation of $L_x$ depends only on the component of that deformation normal to $S_\gamma$. So, the linearization is naturally identified with the map
\[
H^0(L_x,N_{L_x/T})\to H^0(L_x,N_{S_\gamma/T}|_{L_x}).
\]
This map is surjective, as follows from the long exact cohomology sequence associated with the normal bundle exact sequence, since $H^1(L_x, \cO)=0$. Moreover, we already showed the identification $H^0(L_x,N_{S_\gamma/T}|_{L_x})\simeq H^0(\mathbb{CP}^1,\cO(2))$ which is nonzero. Hence the linearization is nonzero (for at least one $u$, but in fact for all transversal to $\gamma$) and \eqref{eq:linear} implies $\xi_x\neq0$, so $c\neq0$. Consequently we compute that $g^{-1}(\xi_x,\xi_x)=c^2g(v,v)\neq0$.
\end{proof}

\subsection*{Proof of Theorem \ref{thm:main}}
The first part follows directly from Theorem \ref{thm:gCKY}. The second one follows from Lemma \ref{lemma:l1}. Indeed, for a conformal geodesic $\gamma$ we get $(\omega,\nabla^E)$ satisfying the gauged CKY equation \eqref{eq:gCKY} as well as \eqref{eq:ASD}. It is sufficient to show that the associated curvature of $\nabla^E$ vanishes along $\gamma$. This is again a consequence of Theorem \ref{thm:gCKY} which characterizes solutions corresponding to conformal geodesics in terms of the curvature.

\end{document}